\documentclass{ifacconf}  
\usepackage{natbib}
\usepackage{algorithm,bbm}
\usepackage[noend]{algpseudocode}
\makeatletter
\def\BState{\State\hskip-\ALG@thistlm}
\makeatother
\usepackage{xcolor}
\usepackage{amssymb,latexsym,amsmath,array,graphicx,mathtools}
\usepackage[utf8]{inputenc}
\usepackage[T1]{fontenc}
\usepackage[ampersand]{easylist}
\usepackage{stfloats,dsfont}
\usepackage{mathrsfs}

\newtheorem{innercustomprop}{Proposition}
\newenvironment{customprop}[1]
  {\renewcommand\theinnercustomprop{#1}\innercustomprop}
  {\endinnercustomprop}
  
  \newtheorem{innercustomlem}{Lemma}
\newenvironment{customlem}[1]
  {\renewcommand\theinnercustomlem{#1}\innercustomlem}
  {\endinnercustomlem}
  
  \newtheorem{innercustomrem}{Remark}
\newenvironment{customrem}[1]
  {\renewcommand\theinnercustomrem{#1}\innercustomrem}
  {\endinnercustomrem}

\newenvironment{proof}{{\em Proof:}}{\hfill$\fbox{}$}
\usepackage{mathtools}

\usepackage{scalefnt} 
\usepackage{blindtext}
\newcommand\scalemath[2]{\scalebox{#1}{\mbox{\ensuremath{\displaystyle #2}}}}
\usepackage[ampersand]{easylist}

\usepackage{mathrsfs}
\DeclareMathAlphabet{\mathpzc}{OT1}{pzc}{m}{it}
\usepackage{subcaption}

\begin{document}
\begin{frontmatter}
\title{
Internal stabilization of an underactuated linear parabolic system via modal decomposition
}

\author[First]{Constantinos Kitsos}    
\author[First]{Emilia Fridman}           

\address[First]{School of Electrical Engineering, Tel-Aviv University, Tel-Aviv, Israel (emails: \{constantinos, emilia\}@tauex.tau.ac.il)  }


\begin{abstract}
This work concerns the internal stabilization of underactuated linear  systems of $m$ heat equations in cascade, where the control is placed internally in the first equation only and the diffusion coefficients are distinct. Combining the modal decomposition method with a recently introduced state-transformation approach for observation problems, a proportional-type stabilizing control is given explicitly. It is based on a transformation for the ODE system corresponding to the comparatively unstable modes into a target one, where calculation of the stabilization law is independent of the arbitrarily large number of them and it is achieved by solving generalized Sylvester equations recursively. This provides a finite-dimensional counterpart of a recently introduced infinite-dimensional one, which led to Lyapunov stabilization. The present approach answers to the problem of stabilization with actuators not appearing in all the states and when boundary control results do not apply.
\end{abstract}
\end{frontmatter}



\section{Introduction}
The problem of control of systems of coupled parabolic PDEs where controls are not placed in all the states, has incited attention since it was posed as an open problem in \cite{Zuazua (2007)}. Such systems have been studied theoretically meanwhile, see survey \cite{Ammar-Khodja et al (2011)}, which collects the plethora of recent studies with respect to various notions of controllability of coupled systems with internal controls. Some answers to  dual to stabilization problems of observer design with less observations than the number of the states have been given in \cite{Kitsos et al (2021a)} and it was seen that the problem of stabilizing $m$ equations with one control only, capturing also the presence of distinct diffusion coefficients, as the one solved here, is quite challenging.

The motivation behind the class of coupled systems that we study comes from biological predator-prey models and other population and social dynamics phenomena in their linearized versions (see \cite{Britton (1986)}). 
Also, in chemical processes (see \cite{Orlov and Dochain (2002)}), coupled temperature-concentration parabolic PDEs occur to describe the process dynamics. Boundary stabilization of parabolic systems with same number of controls as the number of equations has been achieved via backstepping, see \cite{Vazquez and Krstic} and via modal decomposition in \cite{Katz et al (2021)}. For the case of underactuated systems, boundary stabilization has been achieved via backstepping for some classes of hyperbolic systems, see \cite{Coron et al (2013),Aamo (2013),Di Meglio et al (2013)}, however for the class of parabolic systems, strict assumptions should be imposed on the internal dynamics, see \cite{Baccoli et al (2015)}. In the latter, boundary stabilization was achieved only for $2\times 2$ parabolic systems if the second equation is asymptotically stable up to the coupling term. The difficulty in achieving boundary stabilization for general classes of underactuated systems is revealed in controllability studies, see \cite{Ammar-Khodja et al (2011)}. In such works, it is seen that \emph{approximate controllability} leading to boundary stabilization is only possible under specific assumptions, such as when all diffusion coefficients are identical. For the general case of $m$ equations with less than $m$ controlled ones and in the presence of distinct diffusion coefficients, only internal controllability is possible leading to internal stabilization. Thus, to deal with the general case, we can apply stabilizing controls internally. Various studies on internal controllability of underactuated systems (see \cite{Guerrero (2007),Duprez and Lissy (2016),Alabau-Boussouira et al (2017),Fernandez-Cara et al (2015)}) reveal that the difficulty grows with the number of the states and the number of distinct diffusion coefficients as a result of the notion of \emph{algebraic solvability} (see \cite{Steeves et al (2019)}). The problem of stabilization of such systems runs deep, see for instance \cite{Coron (2007),Barbu et al (2006),Barbu (2010),Munteanu (2019)}, see also the seminal work of \cite{Christofides (1998)} on distributed control. However, to the best of our knowledge, explicit stabilizing laws for these underactuated systems have not been determined yet. 

The problem of internal stabilization for the parabolic system in cascade that we study here becomes more difficult for more than $2$ states and distinct diffusion coefficients. The determination of the stabilizing control relies on an existence of a coordinates transformation. In \cite{Kitsos et al (2021a)}, a dual internal observer design problem for parabolic systems was studied for the $3\times 3$ nonlinear case while in \cite{Kitsos et al (2021b)} and \cite{Kitsos (2020)}, the cases of $3\times 3$ and $m\times m$ linear inhomogeneous hyperbolic systems were studied. In these works, appropriate  infinite-dimensional state transformations were introduced to deal with distinct elements on the diagonal of the coefficient of systems' differential operators requiring  the use of higher-order spatial derivatives of the measurement. The reader may also refer to \cite{Alabau-Boussouira (2002)} for the case of underactuated coupled hyperbolic systems and to \cite{Alabau-Boussouira (2013),Alabau-Boussouira (2014),Alabau-Boussouira (2015)} for the case of underactuated hyperbolic and parabolic coupled
cascade systems.

In the present work, to solve similar problems, we follow a modal decomposition approach. We generalize the method mainly used for the scalar case (see \cite{Coron and Trelat (2004)} on direct Lyapunov method for state feedback, see also \cite{Barbu et al (2006)}) or for the vector case but with controls placed in all states (see \cite{Katz et al (2021)} for boundary stabilization), to the case of underactuated systems. We assume that the number of internal inputs is equal to the number of unstable modes and that the resulting matrix that multiplies control inputs in the unstable modes is non-singular. We then introduce a novel state transformation for ODEs with dimension equal to the number of coupled PDEs and written as a polynomial matrix in the slower eigenvalues (needed to be stabilized at given decay rate) of the parabolic operator with order related to the number of distinct diffusion coefficients. The coefficients of this polynomial matrix are nilpotent matrices up to identity matrix, which are subject to recursive generalized Sylvester equations and can be easily determined via a provided algorithm, while their values depend on the dynamics of the parabolic system. The stabilizing law simply relies on the determination of stabilizing gains for the zero-order coupling matrix and also on our introduced state transformation depending on the dynamics only.
In this way, any system specification can lead to a construction of unified and scalable control laws independently of the number of eigenvalues needed to be stabilized, which can be arbitrarily large. 


The paper is organized as follows. The system and conditions along with the description of the problem of stabilization are presented in Section \ref{sec:problem}. The main stabilization strategy is presented in Section \ref{sec:transformation}, where Proposition \ref{proposition} provides our main result. In Section \ref{sec:simul},  we  apply  our  methodology  to  an unstable parabolic system and in Section \ref{conclusion} we provide conclusions and perspectives. 

\textit{Notation:} For a given $x\in {{\mathbb{R}}^{m}}$, $\left| x \right|$ denotes its usual Euclidean norm. For a given matrix $A$ in ${{\mathbb{R}}^{m\times m}}$, ${A}^{\top}$ denotes its transpose, $\left| A \right|:=\sup \left\{ \left| Aw \right|,\left| w \right|=1 \right\}$ is its induced norm, $\mathrm{Sym}(A)=\frac{A+A^{\top}}{2}$ stands for its symmetric part and $\lambda_{\min}(A)$, $\lambda_{\max}(A)$ denote its minimum and maximum eigenvalue, respectively. By $\text{col}\{ a_1,\ldots,a_m \}$ and by $\text{diag}\{ a_1,\ldots,a_m \}$ we denote the column vector and the diagonal matrix, respectively, with elements $a_1,\ldots,a_m.$ By $I_m$ we denote the identity matrix of dimension $m$. For $f, g$ in $L^2\left( 0,L;\mathbb R^m\right)$, by $\left <f,g\right >$ we denote the inner product $\left <f,g \right >=\int_0^L f^\top (x)g(x) dx$ with induced norm $\Vert \cdot \Vert_{L^2\left( 0,L;\mathbb R^m\right)}$, where $L^2\left( 0,L;\mathbb R^m\right)$ denotes the space of equivalence classes of measurable functions $f:[0,L]\to\mathbb{R}^m$. By $\ell^2(\mathbb N;\mathbb R^m)$ we denote the Hilbert space of the square summable sequences $x=(x_n)_{n=1}^{+\infty}$. 
By $\delta_{ij}$  we denote the Kronecker delta $\delta_{ij} = 1$, if $i = j$ and $\delta_{ij} = 0$, otherwise and $\lceil x\rceil $ stands for the ceiling function of $x\in \mathbb R.$
\section{Problem statement}\label{sec:problem}
\subsection{System and requirements}
Consider a system of $m$ coupled parabolic PDEs in a finite domain with controls acting on the first equation only, written as follows $(t,x)\in[0,+\infty)\times(0,L)$:
\begin{align}
\begin{aligned}
    z_t(t,x)=&D z_{xx}(t,x)+Q z + B \sum_{j=1}^Nb_j(x)u_j(t)\\
    z_x(t,0)=&0, \quad \gamma_1 z(t,L)+\gamma_2 z_x(t,L)=0
\end{aligned}\label{sys}
\end{align}
where $z=\text{col}\{z_1, \ldots, z_m\}$ is the state and $D=\text{diag}\left\lbrace d_1,\ldots,d_m\right\rbrace$ is the diffusion matrix
with diffusion coefficients $d_1,\ldots,d_m>0$. The  
 coupling and control coefficients $Q$ and $B$ are assumed to be of the form 
\begin{gather} \label{matrixQ}
Q= 
{
\begin{pmatrix}
q_{1,1}&&\cdots &&q_{1,m}\cr 
q_{2,1}&
 & &&\cr
0&&&&\vdots\cr
\vdots & \ddots &\ddots &\ddots \cr
&&&&\cr
0&\cdots & 0& q_{m,m-1}&q_{m,m}
\end{pmatrix}\nonumber
}, \quad B=\begin{pmatrix}
1\\0\\\vdots\\\vdots\\\vdots\\0
\end{pmatrix},\notag
\end{gather}
meaning that the scalar control actions $u_1(t),\ldots, u_{N}(t)$ with $N$ to be determined later,
act on the first equation only. We assume that the elements of the subdiagonal of $Q$ satisfy
\begin{gather}
q_{2,1}, q_{3,2}, \ldots, q_{m,m-1}\neq 0, \label{controllability}
\end{gather}
which stands as a controllability condition for the pair $(Q,B)$.
Functions $b_1(\cdot),\cdots,b_{N}(\cdot)\in L^2(0,L)$ describe how the control actions are distributed in $[0,L]$ and are subject to some constraints given below. System is subjected to mixed Neumann and Robin conditions, where scalars $\gamma_1, \gamma_2$ satisfy $\gamma_1^2+\gamma_2^2 \neq 0$. 

Before presenting the stabilization method, let us consider the Stürm-Liouville eigenvalue problem
\begin{align}
    d_i \varphi^{\prime \prime}+\bar \lambda \varphi =&0, \quad 0<x<L, \quad i=1,\ldots,m\notag\\
    \varphi^\prime(0)=&0, \quad \gamma_1 \varphi(L)+\gamma_2\varphi^{\prime}(L)=0 \label{SLprob}
\end{align}
admitting a sequence of eigenvalues $\bar \lambda_{n,i}=d_i  \lambda_n$, where $\lambda_n$ are the eigenvalues of \eqref{SLprob} with $d_i=1$. This sequence of eigenvalues corresponds to a sequence of eigenfunctions $(\varphi_n)_{n=1}^{+\infty}$. The eigenvalues form an unbounded increasing and non-negative sequence while the eigenfunctions form a complete orthonormal system in $L^2(0,L)$. Each of the states of \eqref{sys} can be presented as 
\begin{align}
    z_i(t,\cdot)=\sum_{n=1}^\infty z_{i,n}(t) \varphi_n(\cdot), \quad i=1\ldots,m
\end{align}
with coefficients $z_{i,n}$ given by 
\begin{gather}
z_{i,n}=\left <z_i,\varphi_n\right >\label{zi}.
\end{gather}
Taking the time-derivative of \eqref{zi}, substituting dynamics \eqref{sys}, and integrating by parts, we get the following dynamics for $z_n=\text{col}\{z_{1,n},\ldots z_{m,n}\}:$  
\begin{align}
    &\dot z_n(t)=\int_0^L z_t(t,x) \varphi_n(x)  dx\notag\\&=\left (D z_x(L)\varphi_n(L)-Dz(L)  \varphi_n^{\prime} (L)\right )\notag\\&+\left (-\lambda_n D +Q\right )z_n(t)+B \sum_{j=1}^N u_j(t)  \int_0^L \varphi_n(x)  b_j(x) dx, \notag
\end{align}
which by virtue of homogeneous boundary conditions for $\varphi_n(x)$ and $z(t,x)$, is written as follows:
\begin{align}
     \dot z_n(t)=&\left (-\lambda_n D +Q\right )z_n(t) +B \sum_{j=1}^N b_{j,n} u_j(t) \label{zn_dot}
\end{align}
with $b_{j,n}:=\int_0^L   b_j(x) \varphi_n(x) dx$.
Now, given a desired decay rate $\delta>0$, by virtue of the countability and monotonicity of eigenvalues of the parabolic operator,  we can always find a $N \in \mathbb N$  large enough such that 
\begin{align}
-\lambda_{N+1} D +\text{Sym}(Q)+\delta I_m <0.\label{lambdan}
\end{align} 
By monotonicity of $\lambda_n$, the above implies that
\begin{align}
-\lambda_{n} D +\text{Sym}(Q)+\delta I_m <0, \quad \forall n\geq N+1. \label{lambdan2}
\end{align}  
Using the notation $z^N=\text{col}\{z_1,\ldots,z_N\}$ in $\mathbb R^{mN}$, we obtain the following system corresponding to the finite-dimensional part of the eigenspectrum of the parabolic operator:
\begin{align}
\dot z^N(t)=A z^N(t)+\tilde B u(t),  \label{z^Nopen}
\end{align}
where $u(t):=\text{col}\{u_1(t),\ldots,u_N(t)\}$,
\begin{align} 
    A:=&\text{diag}\{-\lambda_1 D+Q,\ldots,-\lambda_N D+Q \}, \label{Amatrix}
\end{align}
and  $\tilde B$ in $\mathbb R^{mN\times N}$ is given by $
\tilde B:= \text{col} \{ B \mathcal B_1,\ldots, B\mathcal B_N\} $
with $\mathcal B_n$ row vectors containing the $n$-th projections of all $b_j(x)$ and written in the form $
\mathcal B_n:=\begin{pmatrix}
b_{1,n}&\cdots&b_{N,n}
\end{pmatrix}, n\geq 1.$
Now, define  matrix 
\begin{align}
\mathcal B_{N\times N}:=&\begin{pmatrix}\mathcal B_1\\\vdots\\\mathcal B_{N}
\end{pmatrix}. \label{mathcalB}
\end{align}
By invoking the Hautus lemma, it is easy to see that the pair $(A, \tilde B)$ is stabilizable whenever the following hypothesis is satisfied on functions $b_j(\cdot)$: 

\textbf{(H)} Matrix $\mathcal B_{N\times N}$ is non-singular. 

Similar hypothesis as (H) appears in several works in the context of stabilization of parabolic systems, see for instance \cite{Hagen and Mezic (2003)}.

We now seek for feedback controls of proportional type as
\begin{align} \label{control}
u_j(t)=&K_j z^N(t),
\end{align}
where $K_{j}\in \mathbb R^{1\times mN}$ are controller gains to be found below.
Then, a direct stabilization approach of system \eqref{z^Nopen} would require to solve inequality
$
\text{Sym} \left (\tilde P(A+\tilde B K) \right)+\delta \tilde P<0, $
where $\tilde P$ in $\mathbb R^{mN\times mN}$ symmetric positive definite and $K:=\text{col}\{K_1,\ldots, K_N\}$. The above is written  in the design LMI form 
\begin{align}
\text{Sym}\left ( AR +\tilde B Y\right )+\delta R<0, \label{LMIfull}
\end{align} 
where we denote the unknowns $R=\tilde P^{-1}$ and $Y=KR$. Then, the desired gain matrix is given by $K=Y R^{-1}.$
This LMI involves matrices of dimension $mN$. 

In this work, we aim at reducing the dimension of the stabilization from $mN$, which depends on the number $N$ of modes to be stabilized, to just the dimension $m$ of the coupled parabolic system, which is fixed. It turns out that this requirement of stabilization is not directly met as a consequence of the presence of distinct diffusion coefficients $d_i$. In fact, we seek for stabilizing actuations $u_j(t)$ whose calculation up to an inversion of matrix $\mathcal B_{N\times N}$ does not depend on the number of the modes $N$ but only on the number of system's equations $m$. Such property is important when dealing with large instabilities in the dynamics or when one would need to efficiently tune the decay rate. In other words, stabilization of \eqref{z^Nopen} should lie on the stabilization of an $m \times m$ matrix, namely, matrix $Q$ and not on each of the diagonal elements of $A$, which can be arbitrarily many depending on the number of modes we need to stabilize at a given rate $\delta$. 
In the next subsection, we will show via examples why stabilization of \eqref{z^Nopen} is not directly implementable when diffusion coefficients are distinct.

At this point, let us   denote 
\begin{align}
K=\mathcal B_{N\times N}^{-1}\text{diag}\{\bar K_1, \ldots, \bar K_N\}\label{Knproperty}
\end{align}
with $\bar K_1, \ldots, \bar K_N$ in $\mathbb R^{1\times m}$ to be determined later.
Closing the loop with control \eqref{control}, $z^N$ satisfies
\begin{align} \label{z^Nsys}
    \dot z^N(t)=\left ( A+F\right )z^N(t),
\end{align}
where $A$ is given by \eqref{Amatrix} and
\begin{align} 
F:=& \text{diag} \{ B \bar K_1,\ldots B \bar K_N\}.\label{A,barb}
\end{align}
This block-diagonal form is appropriate for stabilization as we show in the next sections.

\subsection{The problem of stabilization of the unstable modes} \label{subsec:stabilization}

We present here some scenarios of stabilization of the finite-dimensional part of the eigenspectrum decomposition revealing its difficulty when diffusion matrix $D$ has distinct elements. 
Let us consider \eqref{z^Nsys}. To achieve exponential stability of this system with decay rate $\delta$, one would need to stabilize each of the components $-\lambda_n D+Q$ of matrix $A$ at this rate by choice of appropriate gains $\bar K_n$ as in \eqref{LMIfull}-\eqref{Knproperty}. However, this stabilization strategy would require stabilization of an $mN\times mN$ matrix, which is inefficient when $N$ becomes large.
In order to reduce the stabilization problem for all $N$ modes to just the stabilization of the coupling matrix $Q$, we need to follow an indirect strategy. Indeed, following a direct approach and trying to stabilize only matrix $Q$, one would choose gains
$ \bar K_n=K_Q, n=1,\ldots, N,$ where
$K_Q$ in $\mathbb R^{1\times m}$ is chosen such that $\bar Q:=Q+BK_Q$ satisfies a
Lyapunov matrix inequality of the form  
\begin{align}\label{lyapmatrix}
\text{Sym}\left ( P\bar Q\right )+qP<0
\end{align}
for $P$ in $\mathbb R^{m\times m}$
symmetric positive definite, which is non-diagonal, and $q>0$. This is always possible due to the controllability of $(Q,B)$.
Then, to check asymptotic stability of system \eqref{z^Nsys}, choose Lyapunov function of the form
\begin{align} \label{lyapfunctinitial}
V_0(t)=\frac{1}{2}(z^N(t))^{\top} \bar P z^N(t)
\end{align}
with $ \bar P=\text{diag}\{ P, \ldots, P\}$ consisting of $N$ diagonal blocks $P$.
Then, observe that $\text{Sym} \left (P\left ( -\lambda_n D+\bar Q\right ) \right )$ appearing when taking the time-derivative of $V_0$ is of indefinite sign since $D$ and $P$ do not commute when $D$ has distinct diffusion coefficients and because $P$ is non-diagonal. This means that a stabilizing law chosen to stabilize $Q$ would not automatically lead to the stabilization of all the modes we need to stabilize at rate $\delta$. Note that this complication arising from the lack of a commutative property between the coefficient of the differential operator (the diffusion matrix $D$ here) and a Lyapunov matrix $P$ has been tackled in \cite{Kitsos (2020)}.
To understand how the number of distinct diffusion coefficients plays a role in the complexity of the problem, let us see the following examples.
\subsubsection{Example 1. }\label{exmpl1} Assume that all diffusion coefficients $d_i$ are identical. Then, the stabilization problem would be trivial. Indeed, the gains of the stabilization law \eqref{control} via \eqref{Knproperty} can be chosen  as $\bar K_n=K_Q,$ for all $n=1,\ldots,N$,
where $K_Q\in \mathbb R^{1\times m}$ is chosen such that $\bar Q:=Q+B K_Q$ satisfies \eqref{lyapmatrix} for $P$ symmetric positive definite and $q>0$ sufficiently large depending on the choice of the desired decay rate $\delta$. Then, by choice of Lyapunov function \eqref{lyapfunctinitial}, matrix $\text{Sym} \left (P\left ( -\lambda_n D+\bar Q\right ) \right )=\text{Sym} \left (P\left ( -\lambda_n d_mI_m+\bar Q\right ) \right )$ is always negative definite and the decay rate of system \eqref{z^Nsys} can attain value $\delta$ after appropriate choice of $q$  ($q \geq \delta-\lambda_1 d_m$).

\subsubsection{Example 2. }\label{exmpl2} Let us now see the case where diffusion coefficients are identical up to the first one, namely, $$d_1\neq d_2=\ldots=d_m.$$
We choose gains $\bar K_n$ in \eqref{Knproperty} given by 
\begin{align}
    \bar K_n=-G_n+K_Q, \quad \forall n\in \{1,\ldots,N\}, \label{K_nsigma=2}
\end{align}
where $G_n:=\lambda_n\left (
    d_2-d_1\right )B^\top$.
Again, the gain $K_Q\in \mathbb R^{1\times m}$ is chosen to satisfy \eqref{lyapmatrix} with $\bar Q:=Q+B K_Q$ and by choice of Lyapunov function \eqref{lyapfunctinitial}, system \eqref{z^Nsys} is stabilized at rate $\delta$. This is possible by noting that matrix $\text{Sym} \left (P\left ( -\lambda_n D+Q+B \bar K_n\right ) \right )$, which by \eqref{K_nsigma=2} is equal to $\text{Sym} \left (P\left ( -\lambda_n d_mI_m +\bar Q\right ) \right )$, is negative definite by \eqref{lyapmatrix} and the decay rate of system \eqref{z^Nsys} can be equal to $\delta$ by appropriate choice of $q$ ($q \geq \delta-\lambda_1 d_m$). 
 
\subsubsection{Example 3. }\label{exmpl3}  We finally consider a $3 \times 3$ parabolic system ($m=3$) with $d_2 \neq d_3.$ 
Here, we might have $2$ or $3$ distinct diffusion coefficients and this stabilization problem turns to be more complicated than the previous ones. Indeed, to utilize the previously described Lyapunov stabilization for the finite-dimensional part of the system, we first need to perform a transformation of the form $y_n=T_n z_n, $ for $n=1,\ldots,N$ with 
\begin{align}
    T_n=I_3+\lambda_n \begin{pmatrix}
    0 & \kappa& 0\\0&0&0\\0&0&0
   \end{pmatrix}; \quad \kappa:=\frac{d_2-d_3}{q_{21}}. \label{T_n3}
\end{align}
Then, $y^N=\text{col}\{y_1,\ldots,y_N\}$ in $\mathbb R^{3N}$ satisfies
\begin{align}  
    \dot y^N(t)=\left (\tilde A +\tilde F \right )y^N(t),
\end{align}
where
\begin{align}
   \begin{aligned}
\scalemath{1}{\tilde A:=}&\scalemath{0.9}{\text{diag}\{-\lambda_1 d_3 I_3+Q +BG_1,\ldots, -\lambda_N d_3I_3+Q +BG_N \}; }\\
G_n:=&\begin{pmatrix} \lambda_n \left ( d_3-d_1+\kappa q_{21}\right )\\ \lambda_n^2 \kappa\left ( d_1-d_2-\kappa q_{21}\right )+\lambda_n \kappa \left ( q_{22}-q_{11}\right )\\\lambda_n \kappa q_{23}\end{pmatrix}^\top \notag
\end{aligned}
\end{align}
and $\tilde F:=\text{diag}\{ B \bar K_1 T_n^{-1},\ldots, B \bar K_N T_n^{-1}\}$.
Then, the stabilizing gains are chosen to be of the form
\begin{align*}
\bar K_n=\left ( -G_n+K_Q \right ) T_n,
\end{align*}
where, the first term is needed to eliminate the undesired terms $B G_n$ and, as in the previous examples, $K_Q\in \mathbb R^{1\times m}$ is chosen to satisfy \eqref{lyapmatrix} with $\bar Q:=Q+B K_Q$ and $q$ large enough ($q \geq \delta-\lambda_1 d_3$). 

The above examples demonstrate that the problem of stabilization of an underactuated system is more intricate when diffusion coefficients are distinct, particularly when we have more than two distinct ones. In fact, index
\begin{align}\label{sigma}
\sigma:=&\min\left\lbrace i: d_i = d_j, \forall j=i,i+1,\ldots,m\right\rbrace
\end{align}
is an indicator of the complexity of the problem of stabilization. The larger the value of $\sigma$ is, the more complex is to determine the stabilization law. In our previous examples, for system in Example 1, $\sigma$ was equal to 1 (one diffusion), while in Example 2, $\sigma$ was equal to $2$. Example 3 with $\sigma=3$ provides us with intuition on which indirect strategy we should follow for systems with $\sigma>3$.
In the next section we provide a stabilization law for all possible values of $\sigma$ by determining a state transformation similarly as in the preliminary form \eqref{T_n3}. 

\section{Transformation and stabilizing law}\label{sec:transformation}

In this section, we demonstrate the stabilization of the underactuated system \eqref{sys} via modal decomposition approach. We aim at constructing gains $K_n$ that lead to a closed-loop system, for which we can prove exponential stability. Our idea is to reduce the problem of stabilization for the $mN \times mN$ system to a stabilization problem for system of dimension as large as $m$. We seek for a state transformation that transforms system \eqref{zn_dot} into a target one where this type of control can be applied.

\subsection{Target system and main result}
We first present the target system that allows the derivation of the stabilization law.

 Let us apply a transformation $y_n=T_n z_n$ to system \eqref{zn_dot} with $T_n \in \mathbb R^{m\times m}$ invertible and assumed to be given by
\begin{align} 
T_n=\left\{ \begin{array}{ll} 
I_m+ \sum_{i=1}^{\bar \sigma} \lambda_n^i \bar T_i, \quad & 1\leq  n \leq N,\\I_m, \quad &n\geq N+1\end{array}\right., \label{T_n}
\end{align}
where $$\bar \sigma:=\min\{2\sigma-3,2m-4\}$$ with $\sigma$ given by \eqref{sigma} and $\lambda_n^i$ denoting the $i$-th power of $\lambda_n$. Note that $\mathcal T:=(T_n)_{n=1}^{+\infty}:\ell^2(\mathbb N;\mathbb R^m)\to \ell^2(\mathbb N;\mathbb R^m)$ is a bounded operator with bounded inverse. Matrices $\bar T_i\in \mathbb R^{m\times m}$ are assumed to be of the following nilpotent form for all $i \in \{1,\ldots,\bar \sigma\}:$
\begin{align}
 \scalemath{0.73}{\bar T_i = \begin{pmatrix}
0&\cdots&0&\kappa_{1,\lceil\frac{i}{2}\rceil+1}^i &\kappa_{1,\lceil\frac{i}{2}\rceil+2}^i &\cdots&\cdots&\kappa_{1,m}^i \cr
0&\cdots&0&0&\kappa_{2,\lceil\frac{i}{2}\rceil+2}^i &\cdots&\cdots&\kappa_{2,m}^i\cr
\vdots&&\ddots&&\ddots\cr
0&\cdots&&0&0&\kappa_{m-2-\lceil\frac{i}{2}\rceil,m-2}^i&\kappa_{m-2-\lceil\frac{i}{2}\rceil,m-1}^i&\kappa_{m-2-\lceil\frac{i}{2}\rceil,m}^i\cr
0&\cdots&&0&\cdots&0&\kappa_{m-1-\lceil\frac{i}{2}\rceil,m-1}^i&\kappa_{m-1-\lceil\frac{i}{2}\rceil,m}^i\cr
0&\cdots&&0&\cdots&0&0&0\cr
\vdots&&&\vdots&&\vdots&\vdots&\vdots\cr
0&\cdots&&0&\cdots&0&0&0
\end{pmatrix} }\label{matrixTi}
\end{align}
where $\kappa^i_{j,k}$ are some constants to be determined explicitly in the following, which strictly depend on the dynamics of \eqref{sys} and not on $\lambda_n$. Note that superscripts $i$ appearing in $\kappa^i_{j,k}$ represent indices referring to each of the matrices $\bar T_i$ and not powers while their subscripts $(j,k)$ refer to their position in matrices $\bar T_i$.
By use of this transformation, we aim at obtaining a target system, which after injection of control \eqref{control} and by use of \eqref{Knproperty}, it is written in the closed-loop form
\begin{align}  \left\{ \begin{array}{ll}  \begin{aligned}
    \dot y_n(t)&=\left (-\lambda_n d_m I_m +Q +B G_n \right.\\&\left.+ B \bar K_nT_n^{-1}\right )y_n(t),  \qquad n \leq N,\\
    \dot y_n(t)&=\left (-\lambda_n D+Q\right )y_n(t) \\&+B \sum_{j=1}^N b_{jn} K_j z^N(t),  \quad  n\geq N+1
    \end{aligned}\end{array}\right.  \label{target_sys} 
\end{align}
with 
$G_n$ given by 
\begin{align}
 &G_n=-B^\top \left(\left ( Q-\lambda_n d_m I_m\right )\left( \sum_{i=1}^{\bar \sigma} \bar T_i \lambda_n^i\right )\right.  \notag\\&\left. +\left( \sum_{i=1}^{\bar \sigma} \bar T_i \lambda_n^i\right )\left ( \lambda_n D-Q\right )+\left ( D-d_mI_m\right )\lambda_n \right)T_n^{-1}.  \label{Gn}
\end{align}
Terms $B G_n y_n(t)$ are undesired in the stabilization process but they can be canceled by use of the gains $\bar K_n$. To obtain target system \eqref{target_sys}, let us assume that $\bar T_i$ satisfy the following recursive generalized Sylvester equations
for all $ i\in \{1,\ldots,\bar \sigma\}$:
\begin{align} 
    \left(  I_m-BB^\top\right ) \left ( Q \bar T_i-\bar T_i Q+ \bar T_{i-1}\left ( D-d_m I_m \right )\right )=0, \label{Sylvester}
\end{align}
where $\bar T_0:=I_m$. We obtain the following result on solutions to \eqref{Sylvester}:
\begin{customlem}{1} \label{lemma}
If condition \eqref{controllability} holds true, there exist matrices $\bar T_i$ of the form \eqref{matrixTi} satisfying generalized Sylvester equations \eqref{Sylvester}. Their components $\kappa^i_{j,k}$ can be obtained explicitly.
\end{customlem}
\begin{proof}
Solvability of this family of generalized Sylvester equations \eqref{Sylvester} is hard to be checked by mathematical criteria, however, we are in a position to directly determine solutions.  Thanks to the special structure of $\bar T_i$, we can apply an elimination procedure of each element of matrix $\left(  I_m-BB^\top\right ) \left ( Q \bar T_i-\bar T_i Q+ \bar T_{i-1}\left ( D-d_m I_m \right )\right )$ in a recursive manner. For each row, we start from elimination of its rightmost element and then we eliminate one by one all of its elements by moving one position to the left. The procedure initiates at the lowest row with nonzero elements and when all elements of the current row are eliminated leftwards, we recede to the rightmost element of one row before it and we continue the same procedure until all elements are eliminated. For each of these eliminations, we calculate an element $\kappa^i_{j,k}$ as the sole unknown in this entry, which is written as a linear combination of elements that have been already calculated in precedent eliminations. One can easily check that extracting a sole unknown component $\kappa^i_{j,k}$ for each of these eliminations is a result of the  special structure of \eqref{matrixTi} and controllability condition \eqref{controllability}. More precisely, Algorithm \ref{algorithm} below describes the procedure to calculate each of the elements $\kappa^i_{j,k}$ of $\bar T_i$.
\begin{algorithm}
\caption{Calculation of transformation $T_n$}\label{algorithm}
\begin{algorithmic}[1]
\Procedure{Compute $\kappa_{j,k}^i,$ for all  $i \in \{1,\ldots,\bar \sigma\}, j\in \{1,\ldots,m-1-\lceil\frac{i}{2}\rceil\},k \in \{\lceil\frac{i}{2}\rceil+1,\ldots,m\}.$}{}
\State {$\bar T_0:=I_m$ and matrices $\bar T_i$ have the form \eqref{matrixTi} for $i\in \{1,\ldots,\bar \sigma\}.$ }
\State{$i=1$.}
\While{$i\leq \bar\sigma,$}
\State{$j=m-1-\lceil\frac{i}{2}\rceil$.}
\While{$j \geq 1$}
\State{$k=m$}
\While{$k\geq j+\lceil\frac{i}{2}\rceil$}
\State{Calculate $\kappa^i_{j,k}$ by eliminating element $(j,k)$ of matrix $\left(  I_m-BB^\top\right ) \left ( Q \bar T_i-\bar T_i Q+ \bar T_{i-1}\left ( D-d_m I_m \right )\right )$. In each step, all calculated $\kappa^i_{j,k}$ are written as linear combinations of $\kappa^i_{j,k}$ already calculated in previous steps.}
\State {$k\gets k-1$.}
\EndWhile
\State {$j\gets j-1$.}
\EndWhile
\State {$i\gets i+1$.}
\EndWhile
\EndProcedure
\end{algorithmic}
\end{algorithm}
By applying this algorithm, we achieve to calculate all $\kappa^i_{j,k}$ appearing in \eqref{matrixTi}. The following formula retrieves each of these elements in the recurring order Algorithm \ref{algorithm} suggests. More precisely, we obtain for all $i \in \{1,\ldots,\bar \sigma\}, j\in \{1,\ldots,m-1-\lceil\frac{i}{2}\rceil\},k \in \{\lceil\frac{i}{2}\rceil+1,\ldots,m\}$
\begin{align}
    &\kappa^i_{j,k}=\frac{1}{q_{j+1,j}}\left ( \sum_{l=0}^{m-j-1-\frac{\lceil{i}\rceil}{2}} \kappa^i_{j+1, j+\frac{\lceil{i}\rceil}{2}+1+l}q_{j+\frac{\lceil{i}\rceil}{2}+1+l,k}\right.\notag\\&\left.-\sum_{l=0}^{m-j-1-\frac{\lceil{i}\rceil}{2}}q_{j+1,j+1+l}\kappa^i_{j+1+l,k}+\kappa^{i-1}_{j+1,k}(d_m-d_{k})\right ), \label{kappaformula}
\end{align}
where we define $\kappa^0_{j+1,k}:=\delta_{j+1,k},$ for all $j\in \{1,\ldots,m-2\}, k\in \{1,\ldots,m\}$. Note here that formula \eqref{kappaformula} is well-defined thanks to the controllability condition \eqref{controllability}.
\end{proof}

We are now in a position to state and establish our main result. 
\begin{customprop}{1}\label{proposition}
Consider parabolic system \eqref{sys} with initial condition $z(0,\cdot)=:z^0(\cdot)$ in $H^1\left (0,L;\mathbb R^m\right)$ satisfying compatibility conditions and assume that controllability condition \eqref{controllability} and hypothesis (H) hold true. Consider also matrices $T_n$ given by \eqref{T_n} with $\bar T_i$ solving \eqref{Sylvester}.  Given a decay rate $\delta>0$, let $N\in\mathbb N$ be subject to \eqref{lambdan}. Assume that there exist $0<W\in  \mathbb R^{m \times m}$ and $Z\in \mathbb R^{1\times m}$ satisfying the following LMI:
\begin{align} \label{omega,psi}
 \text{Sym}\left (QW+BZ\right )+\left(\delta-\lambda_1 d_m \right ) W <0.
\end{align} 
Denote $K_Q=ZW^{-1}$. Let $\bar K_n\in \mathbb R^{1\times m}$ be given by
$\bar K_n=\left (-G_n+K_Q \right)T_n,$
for all $n=1,\ldots, N,$ where via expression \eqref{Gn} for $G_n$, the above is rewritten as
\begin{align}
\bar K_n=B^\top \Big ( \left ( Q-\lambda_n d_m I_m\right ) T_n+T_n\left ( \lambda_n D-Q\right ) \Big ) +K_Q T_n.\label{gains0}
\end{align}
Then, the controller \eqref{control} with gains $K_n$, defined by \eqref{Knproperty}, exponentially stabilizes \eqref{sys} with a decay rate $\delta$, meaning that the solutions of the closed-loop system satisfy the following inequality: 
\begin{align}
    \Vert z(t,\cdot)\Vert_{L^2\left( 0,L;\mathbb R^m\right)} \leq M e^{-\delta t} \Vert z^0(\cdot) \Vert_{L^2\left( 0,L;\mathbb R^m\right)}, \forall t \geq 0 \label{stability0}
\end{align}
with  $M>0$.  Moreover, \eqref{omega,psi} is always feasible. 
\end{customprop} 
The proof of Proposition \ref{proposition} follows in the next subsection. This result illustrates that stabilization relies on the determination of a stabilizing gain $K_Q$ for coupling matrix $Q$  and also on the calculation of a family of matrices $\bar T_i$, whose number depends on the number of distinct diffusion coefficients $d_i$ (represented by $\bar \sigma$) while their values only depend on system's dynamics. These matrices $\bar T_i$ are calculated easily by following the algorithm Lemma \ref{lemma} suggests. This stabilization method is scalable up to the inversion of matrix $\mathcal B_{N\times N}$ given by \eqref{mathcalB}.  Note that \eqref{omega,psi} is the standard design LMI form of Lyapunov inequality \eqref{lyapmatrix} and it is easily solvable by LMI toolboxes.

\subsection{Proof of Proposition \ref{proposition}}
We prove here Proposition \ref{proposition}. 

Note first that transformation $T_n$ appearing in stabilization law is calculated via the constructive Algorithm \ref{algorithm} coming from Lemma \ref{lemma}. To see how $T_n$ maps \eqref{zn_dot} to target system \eqref{target_sys}  via control \eqref{control}, let us apply it to \eqref{zn_dot}. Then, we obtain
\begin{align}
    \dot y_n(t)=&\left (-\lambda_n T_n D T_n^{-1 } +T_n Q T_n^{-1}+B K_n T_n^{-1}\right )y_n(t),
\end{align}
for all $n=1,\ldots,N.$
Comparing the above system with target system \eqref{target_sys}, the following equations must be satisfied for all $n=1,\ldots,N$:
\begin{align}
    \left ( Q-\lambda_n d_m I_m\right )T_n +T_n\left ( \lambda_n D-Q\right )+B G_n T_n=0. \notag
\end{align}
Substituting  \eqref{T_n} in the previous equation, this is written as
\begin{align}
   &\left ( Q-\lambda_n d_m I_m\right )\left( \sum_{i=1}^{\bar \sigma} \bar T_i \lambda_n^i\right )+\left( \sum_{i=1}^{\bar \sigma} \bar T_i \lambda_n^i\right )\left ( \lambda_n D-Q\right )\notag\\&+\left ( D-d_mI_m\right )\lambda_n+B G_n T_n =0, \quad n=1,\ldots,N. \label{sylvesterwithGn}
\end{align}
After injecting expression for $G_n$, \eqref{sylvesterwithGn} is written as 
\begin{align}
    &\left( I_m-BB^\top\right ) \left ( \left ( Q-\lambda_n d_m I_m\right )\left( \sum_{i=1}^{\bar \sigma} \bar T_i \lambda_n^i\right )\right.\notag\\&\left. +\left( \sum_{i=1}^{\bar \sigma} \bar T_i \lambda_n^i\right )\left ( \lambda_n D-Q\right )+\left ( D-d_mI_m\right )\lambda_n\right )=0. \label{firstSYlvester}
\end{align}
Then, eliminating all the coefficients of $\lambda_n^i$ in \eqref{firstSYlvester} for all $i$ in $\{1,\ldots,\bar \sigma\}$, we obtain \eqref{Sylvester}, which is assumed to hold true for all $\bar T_i$, $i\in \{1,\ldots,\bar \sigma\}$. Therefore, \eqref{Sylvester} guarantees that, via transformation $T_n$, we obtain target system \eqref{target_sys}.

At this point, let us remark that for given initial condition $z^0$ in $H^1\left (0,L;\mathbb R^m\right)$ satisfying compatibility conditions for \eqref{sys}, existence of unique classical solutions to system \eqref{sys} with nonlocal feedback control \eqref{control}, namely $z\in C^1\left ([0,+\infty);L^2\left (0,L;\mathbb R^m\right) \right )$, follows from simple argument such as the Lumer-Philipps theorem, see for example (\cite{Pazy (1983)}, Corollary 4.4, Chapter 1.  

Let us now prove $L^2$ stability of the closed-loop system applying direct Lyapunov method (see for example \cite{Katz and Fridman (2020)}). First, observe that by injecting gains \eqref{gains0}, $y^N:=\text{col}\{y_1,\ldots,y_N\}\in \mathbb R^{mN}$ satisfies \begin{align} 
    \dot y^N(t)=H y^N(t), \label{y^Nclosed}
\end{align}
where $H:=\text{diag}\{H_1,\ldots,H_N\}$ with $H_n:=-\lambda_n d_mI_m+ Q+B K_Q, n=1,\ldots,N.$
Now, by the fact that $(Q,B)$ is controllable, we can stabilize matrix $Q$, in such a way that we can find $0<P$ in $\mathbb R^{m \times m}$, $K_Q$ in $\mathbb R^{1\times m}$, and $\bar \rho>0$ such that
\begin{align}
\scalemath{0.95}{\begin{pmatrix} \text{Sym} \left ( P(Q+B K_Q) \right )+\left(\delta-\lambda_1 d_m \right ) P &I_m\\I_m&-\bar \rho I_m \end{pmatrix}<0, }\label{LMIP}
\end{align}
which is written in the design LMI form
\begin{align} 
\scalemath{1}{\begin{pmatrix} \text{Sym}\left (QW+BZ\right )+\left(\delta-\lambda_1 d_m \right ) W &I_m\\I_m&-\bar \rho I_m\end{pmatrix}<0,} \notag
\end{align} 
where $W=P^{-1}$ and $Z=K_Q W$. The latter implies also feasibility of LMI \eqref{omega,psi}.
Next, by virtue of \eqref{lambdan}, we can always find $\rho>0$ such that the following LMI is satisfied:
\begin{align}
\begin{pmatrix}
-\lambda_{N+1} D +\text{Sym}(Q)+\delta I_m & \frac{1}{\sqrt 2}I_m\\\frac{1}{\sqrt 2}I_m&-\rho I_m
\end{pmatrix}<0.  \label{Omega}
\end{align}
To prove stability, defining first $y(t,x)=\sum_{n=1}^{+\infty}\varphi_n(x)y_n(t)$, consider Lyapunov functional $\mathcal V:L^2\left (0, L;\mathbb R^m\right)\to \mathbb R$
\begin{align}
\mathcal V[y]=\frac{1}{2} \left (y^N\right )^\top \bar P y^N+\frac{\rho_0}{2} \sum_{n=N+1}^{+\infty} \vert y_n\vert^2,
\end{align} \begin{align}\label{rho0}\text{where }\rho_0:=\frac{2}{\rho \bar \rho\beta \vert K \vert^2};  \beta:=\max_{n=1,\ldots,N} \vert T_n^{-1}\vert^2 \sum_{j=1}^{N}\Vert b_j(\cdot) \Vert_{L^2 \left (0,L\right )}^2\end{align} with $\rho>0$ satisfying \eqref{Omega}, $\bar \rho>0$ satisfying \eqref{LMIP} and $K$ given by \eqref{Knproperty} (where $\bar K_n$ are given by \eqref{gains0}). Also, $\bar P:=\text{diag}\{P,\ldots,P\}$.  By invoking boundedness of $\mathcal T, \mathcal T^{-1}$ in $\ell^2\left (\mathbb N;\mathbb R^{m}\right)$, the fact that $y_n=T_n z_n,$ the Cauchy-Schwarz inequality, and Parseval's identity, we get
\begin{align*}
\underline c \Vert z(t,\cdot)\Vert^2_{L^2 \left (0,L;\mathbb R^m\right )}=\underline c \sum_{n=1}^{+\infty}\vert z_n(t)\vert^2\leq &\sum_{n=1}^{+\infty} \vert y_n(t)\vert^2 \notag\\&\leq  \bar c \Vert z(t,\cdot)\Vert^2_{L^2 \left (0,L;\mathbb R^m\right )},
\end{align*}
where $\underline c:=\frac{1}{\max_{n\in \mathbb N}\vert T_n^{-1}\vert^2}$ and $\bar c:=\max_{n \in \mathbb N}\vert T_n \vert^2.$
By continuous differentiability of solutions with respect to $t$ for all $t\geq 0$, we are in a position to define $V(t):=\mathcal V[y](t)$ for all $t\geq 0$ and we may take its time-derivative $\dot V(t)$ along the solutions of target system \eqref{target_sys}. By use of the previous inequality, we obtain for $V(t)$
\begin{align}
\frac{\underline c}{2}\min \left (\lambda_{\min}(P),\rho_0\right )\Vert z(t,\cdot) \Vert_{L^2 \left (0,L;\mathbb R^m\right )}^2\leq V(t) \notag\\\leq \frac{\bar c}{2}\max \left (\lambda_{\max}(P),\rho_0\right )\Vert z(t,\cdot) \Vert_{L^2 \left (0,L;\mathbb R^m\right )}^2. \label{Vineq}
\end{align}
Its derivative satisfies
\begin{align}
\dot V(t) =& \left (y^N(t)\right )^\top \text{Sym}(\bar P H)y^N(t)\notag\\&+\rho_0\sum_{n=N+1}^{+\infty}y_n^\top(t) \left( -\lambda_n D + \text{Sym}\left ( Q\right )\right) y_n(t)\notag\\&+\rho_0 \sum_{n=N+1}^{+\infty}y_n^\top(t) B \sum_{j=1}^N b_{jn} K_j z^N(t). \label{Lyapder}
\end{align}
By the Cauchy-Schwarz inequality and Parseval's identity, last term of \eqref{Lyapder} is bounded as follows:
\begin{align*}
&\rho_0 \sum_{n=N+1}^{+\infty}y_n^\top(t) B \sum_{j=1}^N b_{jn} K_j z^N(t) \leq \rho_0 \frac{1}{2\rho}  \sum_{n=N+1}^{+\infty} \vert y_n(t)\vert^2 \notag\\&+\rho_0\frac{\rho}{2}  \vert z^N(t)\vert^2 \vert K\vert^2 \sum_{n=N+1}^{+\infty}\vert \mathcal B_n \vert^2\notag\\&\leq \rho_0\frac{1}{2\rho} \sum_{n=N+1}^{+\infty} \vert y_n(t)\vert^2 +\rho_0 \frac{\rho}{2}\beta \vert K\vert^2 \vert y^N(t)\vert^2
\end{align*}
where $\rho>0$ satisfies \eqref{Omega}  and $\beta$ is given by \eqref{rho0}. After substituting expression \eqref{rho0} for $\rho_0$, \eqref{Lyapder} is bounded as
\begin{align}
\dot V(t)\leq -2 \delta V(t) +\left (y^N(t)\right )^\top \Gamma y^N(t)\notag\\+\rho_0\sum_{n=N+1}^{+\infty} y_n^\top(t) \Omega_n y_n(t), \label{Lyapder2}
\end{align} 
where $\Gamma:=\text{diag}\{\text{Sym}(P H_1)+\left ( \frac{1}{\bar \rho}+\delta P \right )I_{m},\ldots,$\\$\text{Sym}(P H_N)+\left ( \frac{1}{\bar \rho}+\delta P \right )I_{m}\}$ and $\Omega_n:=-\lambda_n D+\text{Sym}(Q)+\left (\frac{1}{2\rho}+\delta \right )I_m$.
Monotonicity of the eigenvalues, in conjunction with \eqref{LMIP} and \eqref{Omega}, implies $\Gamma<0$ and $ \Omega_n<0, n \geq N+1,$ respectively. Thus, \eqref{Lyapder2} in conjunction with \eqref{Vineq} readily yields to a stability inequality of the form \eqref{stability0}. The proof is complete. 

\begin{customrem}{1}
The finite-dimensional transformation $T_n$ \eqref{T_n} is directly related to an infinite-dimensional one firstly introduced in \cite{Kitsos (2020),Kitsos et al (2021a)} for an observer design problem with only one observation corresponding to various classes of coupled PDEs. In these works, the corresponding transformation was a matrix operator with high-order differentiations in its domain. Note that the role of exponents of $\lambda_n$ in \eqref{T_n} were played there by differential operators of the same order as these exponents and the matrix equations \eqref{firstSYlvester} were corresponding to general operator equations. Note also that those works captured space-varying and nonlinear dynamics. Those cases, being more general than the ones here, required strong regularity assumptions and cannot be tackled by the present approach. 
\end{customrem}

\section{Simulation}\label{sec:simul}
\begin{figure}
\includegraphics[scale=0.3]{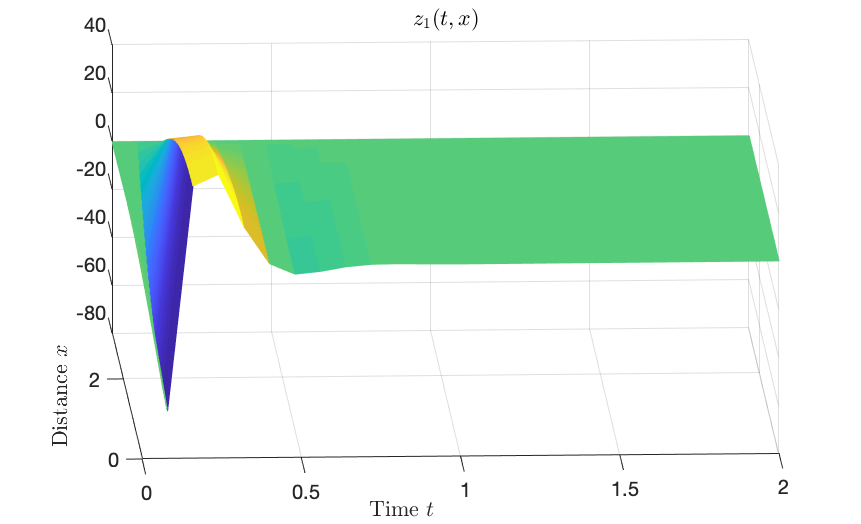}
\caption{Time and space evolution of first state}\label{fig1}
\end{figure}

\begin{figure}
\includegraphics[scale=0.3]{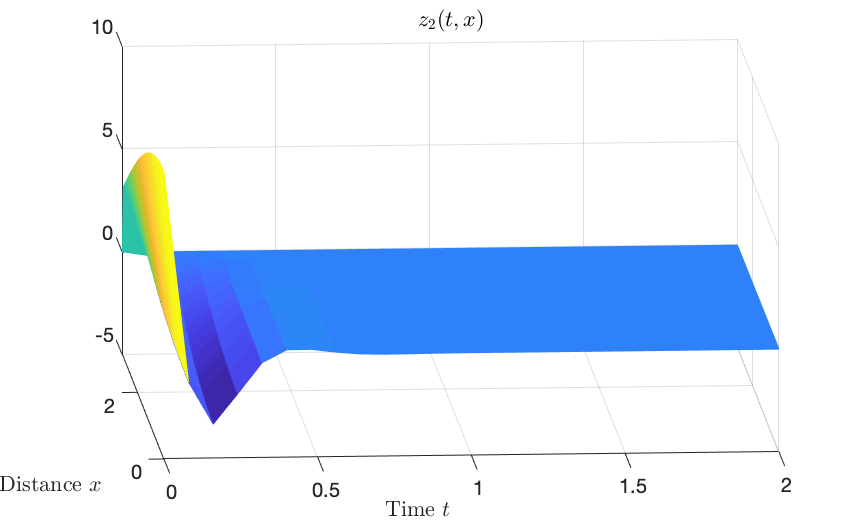}
\caption{Time and space evolution of second state }\label{fig2}
\end{figure}

\begin{figure}
\includegraphics[scale=0.3]{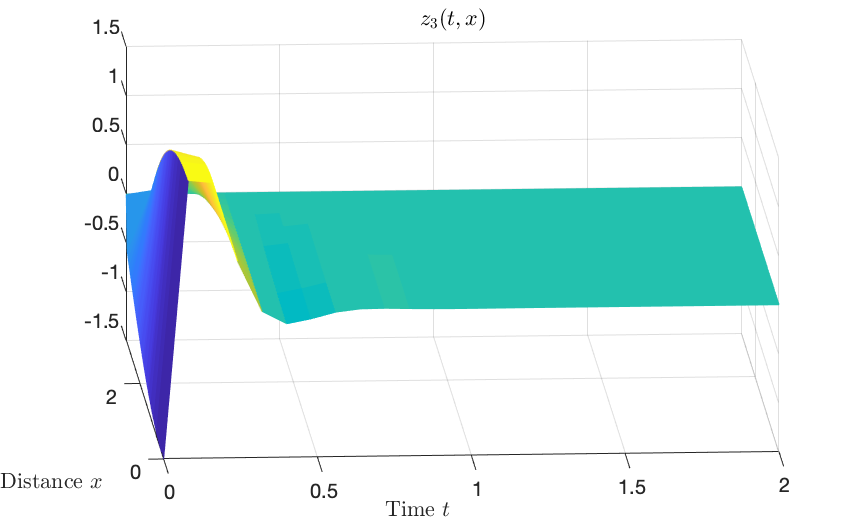}
\caption{Time and space evolution of third state }\label{fig3}
\end{figure}

We consider the stabilization of  an {unstable} parabolic system \eqref{sys} of $3$ equations with $L=\pi$, $\gamma_1=1, \gamma_2=0,$ and $
D=\text{diag}\{4,5,6 \}, $ $ Q=\begin{pmatrix}10&4&8\\ 1&10&2\\0&1&20\end{pmatrix}.
$
We retrieve from the Stürm-Liouville problem \eqref{SLprob} 
$\lambda_n=\left (n-\frac{1}{2}\right )^2 \pi^2/L^2, \varphi_n(x)=\sqrt 2 \cos(\sqrt{\lambda_n}x)$. We choose decay rate $\delta=9$.  We select $N=3$ satisfying inequality \eqref{lambdan}.  Since $m=3$, transformation \eqref{T_n} is of the form \eqref{T_n3} for all $n\leq N$. Shape functions are chosen as $b_j(x)=\mathds 1_{[0.1 j, 0.1 j+0.1]}, j=1,\ldots, N,$ in such a way that matrix $\mathcal B_{N\times N}$ satisfies hypothesis (H) for $N=3$. {We then calculate $K_Q$ by solving LMI \eqref{omega,psi}, which is given by $K_Q=\begin{pmatrix}-67.5&	-3059&-5823\end{pmatrix}$. By applying \eqref{gains0}, we obtain  $\bar K_1=\begin{pmatrix}-68.5	& -2990 &-5821 \end{pmatrix}, \bar K_2=\begin{pmatrix}-59.5	&-1915	-5791 \end{pmatrix}, \bar K_3=\begin{pmatrix}-4.5&	3138&	&-5661 \end{pmatrix}.$
We then calculate the controller gains $K_j$ by using \eqref{Knproperty}. } Simulations of all three PDE states of the closed-loop system are shown in Figures \ref{fig1}-\ref{fig3} for choice of initial condition $z^0(x)=\left(\begin{smallmatrix}\cos x+1\\ 6\cos\frac{x}{2}+3\\-\cos \frac{x}{2}-0.5\end{smallmatrix}\right)$.
Furthermore, by using standard LMI solvers, our method illustrated in Proposition \ref{proposition} to calculate stabilization gains $K_j$ in \eqref{control}  is compared with LMI solving resulting from the direct approach \eqref{LMIfull}. For $N=3$, our indirect approach is approximately $2$ times faster with respect to elapsed time, while for $N=10$ (corresponding to larger $\delta$), it was $270$ times faster than standard LMI.  Note also that for large values of $N$, \eqref{LMIfull} turns out to be computationally hard, while our algorithm to calculate controller gains $K_j$ does not suffer from such limitations. 

\section{Conclusion} \label{conclusion}
The problem of internal stabilization of an underactuated parabolic system in a cascade form and in the presence of distinct diffusion coefficients was considered based on modal decomposition. The stabilization problem was reduced to just the stabilization of the coupling matrix of the parabolic PDEs avoiding in that way a direct stabilization of the whole system of ODEs corresponding to the comparatively unstable modes, which might have arbitrarily large dimension. An easily calculable state transformation of dimension equal to the number of coupled PDEs was introduced in order to solve this problem. 

In our future works, the present approach will be applied to various problems such as non-diagonal cascade systems, higher spatial dimensions, or when the dimension of the finite-dimensional part of the eigenspectrum decomposition of the differential operator plays a crucial role, such as observer-based control, systems with large uncertainties with known bounds, and nonlinear systems.

\end{document}